\documentclass{amsart}
\usepackage{amsmath, amssymb,verbatim}
\newtheorem{theorem}{Theorem}

\newtheorem{lemma}{Lemma}

\author[BOUHALI]{Aissa Bouhali}
\address{Département de Mathématiques, Ecole Normale Supérieure de Laghouat, Algeria.
    Laboratory of Pure and Applied Mathematics, University of Laghouat, Algeria.}
\email{aissa.bouhali@ens-lagh.dz}
\author[LOUHICHI]{Issam Louhichi}
\address{Department of Mathematics \& Statistics, College of Arts \& Sciences, American University of Sharjah, P.O.Box 2666, Sharjah, UAE.}
\email{ilouhichi@aus.edu}
\begin{document}
\title[Commutant of sum of two Toeplitz operators ]
{Commutant of sum of two quasihomogeneous Toeplitz operators}

\date{\today} 
\begin{abstract}
	A major open question in the theory of Toeplitz operator on the Bergman space of the unit disk of the complex plane is to fully characterize the set of all Toeplitz operators that commute with a given one. In \cite{al}, the second author described the sum $S=T_{e^{im\theta}f}+ T_{e^{il\theta}g}$, where $f$ and $g$ are radial functions, that commutes with the sum $T=T_{e^{ip\theta}r^{(2M+1)p}}+ T_{e^{is\theta}r^{(2N+1)s}}$. It is proved that $S=cT$, where $c$ is a constant. In this article, we shall replace $r^{(2M+1)p}$ and $r^{(2N+1)s}$ by $r^n$ and $r^d$ respectively, with $n$ and $d$ in $\mathbb{N}$, and we shall show that the same result holds.
\end{abstract}	
\maketitle
\section{Introduction}
In the complex place $\mathbb{C}$, we denote by $\mathbb{D}$ the open unit disk and $dA=rdr\frac{d\theta}{\pi}$, where $(r,\theta)$ are the polar coordinates, the normalized Lebesgue measure on $\mathbb{D}$. Let $L^2(\mathbb{D},dA)$ be the Hilbert space of all square integrable functions on $\mathbb{D}$ with respect to the measure $dA$. 

The classical unweighted Bergman space $L^2_a(\mathbb{D})$ is the closed subspace of $L^2(\mathbb{D},dA)$ consisting of all analytic functions on $D$. Moreover, the set $\{z^n: n=0,1,2,\ldots\}$ is an orthogonal basis for $L^2_a(\mathbb{D})$. Since $L^2_a(\mathbb{D})$ is closed, the orthogonal projection $P$ from $L^2(\mathbb{D},dA)$ onto $L^2_a(\mathbb{D})$, called the Bergman projection, is well defined. See \cite{h} for more information about the theory of Bergman spaces.

For a function $f\in L^2(\mathbb{D},dA)$, we define the Toeplitz operator $T_f$ from $L^2_a(\mathbb{D})$ into itself by $T_f(g)=P(fg)$ whenever $fg$ is in $L^2_a(\mathbb{D})$. The function $f$ is called the symbol of the Toeplitz operator $T_f$. From our definition of $T_f$, it is clear that bounded analytic functions are in the domain of $T_f$ and so $T_f$ is densely defined on $L^2_a(\mathbb{D})$. Moreover, it is easy to see that if the symbol $f$ is bounded on $\mathbb{D}$, then $T_f$ is bounded and $||T_f||\leq ||f||_{\infty}$.

Over the past half century or so, various algebraic properties of Toeplitz operators have been subject to scrutiny investigations. Nevertheless, the problem of describing the commutant of a given Toeplitz operator, i.e., the set of all Toeplitz operators that commute (in a sense of composition) with it, is still far from being totally solved. Almost nothing is known to us about when $T_fT_g=T_gT_f$ for "general symbols" $f$ and $g$. See \cite{sl, al, acr, cr, lt, l, lr, lry, lz, rv, sz} to know more about the commuting problem of Toeplitz operators on $L^2_a(D)$ and the answers obtained so far for certain classes of symbols.

In this work, we will be dealing with the class of the so-called quasihomogeneous Toeplitz operators. A symbol $f$ is said to be quasihomogeneous of degree an integer $p$ if $f(re^{i\theta})=e^{ip\theta}\phi(r)$, where $\phi$ is a radial function. In this case, the associated Toeplitz operator $T_f$ is also called quasihomogeneous Toeplitz operators of degree $p$ (see \cite{cr, lr, lsz}). The motivation behind the study of such family of symbols is that $L^2(\mathbb{D},dA)$ can be written as $L^2(\mathbb{D},dA)=\bigoplus_{k\in\mathbb{Z}}e^{ik\theta}\mathcal{R}$, where $\mathcal{R}$ is the space of square integrable functions on $\mathbb{D}$ with respect to the measure $dA$. Thus every function $f\in L^2(\mathbb{D},dA)$ has the following polar decomposition $f(z)=f(re^{i\theta})=\sum_{k\in\mathbb{Z}}e^{ik\theta}f_k(r)$, with the $f_k$'s being radial functions. We believe that the study of quasihomogeneous Toeplitz operators would help us to characterize the commutant of Toeplitz operators with more general symbols.

This present work is motivated by \cite[Theorem 3.1,~p.52]{al}. For our purposes we choose to state the latter and it can be summarized as follows:
\begin{theorem}\label{thm4}
	Let $\phi(r)=r^{(2M+1)p}$ and $\psi(r)=r^{(2N+1)s}$ with $p<s$, $M$, and $N$ being all integers greater or equal to $1$. Assume there exist $m, l\in\mathbb{N}$ and nontrivial radial functions $f, g$ such that  the following two hypotheses are satisfied
	\begin{displaymath}
		\left\{\begin{array}{ll}
			(H1)& T_{e^{im\theta}f}+ T_{e^{il\theta}g}
			\textrm{ commutes with } T_{e^{ip\theta}\phi}+T_{e^{is\theta}\psi},\\
			(H2)& 1\leq p<s, \ 1\leq m<l, \textrm{ and } l+p=m+s.
		\end{array}\right.
	\end{displaymath}
	Then $m=p,\ l=s$ and
	$$T_{e^{im\theta}f}+T_{e^{il\theta}g}=c\left(T_{e^{ip\theta}\phi}+T_{e^{is\theta}\psi}\right)$$ for some constant $c$.
\end{theorem}
Our goal here is to replace the radial parts of the symbols $e^{ip\theta}\phi$ and $e^{is\theta}\psi$ in Theorem 1, namely $\phi(r)=r^{(2M+1)p}$ and $\psi(r)=r^{(2N+1)s}$, by $\phi(r)=r^n$ and $\psi(r)=r^d$ with $n,d\in\mathbb{N}$. More precisely, we will get rid off the condition that the power of $r$ in $\phi$ and $\psi$ is an odd number times the quasihomogeneous degrees $p$ and $s$ respectively. In doing so, we will still be able to use the result of the second authors and N. V. Rao in \cite{lrp} about the existence of roots for Toeplitz operators. Their result states that there exist radial functions $\tilde{\phi}$ and $\tilde{\psi}$ such that $T_{e^{ip\theta}r^n}=\left(T_{e^{i\theta}\tilde{\phi}}\right)^p$ and $T_{e^{is\theta}r^d}=\left(T_{e^{i\theta}\tilde{\psi}}\right)^s$. Now if the hypotheses $(H1)$ and $(H2)$ are satisfied, then \cite[Remark 2]{lr} implies that for every vector $z^k$ of the orthogonal basis of $L^2_a(\mathbb{D})$ we have
\begin{eqnarray}
	T_{e^{im\theta}f}T_{e^{ip\theta}r^n}(z^k)&=&T_{e^{ip\theta}r^n}T_{e^{im\theta}f}(z^k),\\
	T_{e^{il\theta}g}T_{e^{is\theta}r^d}(z^k)&=&T_{e^{is\theta}r^d}T_{e^{il\theta}g}(z^k),\\
	\left(T_{e^{im\theta}f}T_{e^{is\theta}r^d}+T_{e^{il\theta}g}T_{e^{ip\theta}r^n}\right)(z^k)
	&=&\left(T_{e^{is\theta}r^d}T_{e^{im\theta}f}+T_{e^{ip\theta}r^n}T_{e^{il\theta}g}\right)(z^k).
\end{eqnarray}
From equations $(1)$ and $(2)$ we can say that $T_{e^{im\theta}f}$ and $T_{e^{il\theta}g}$ commute with $T_{e^{ip\theta}r^n}$ and  $T_{e^{is\theta}r^d}$ respectively. Thus \cite[Proposition 2 and Lemma 2]{lr} imply
\begin{equation}\label{m}
T_{e^{im\theta}f}=c_1\left(T_{e^{i\theta}\tilde{\phi}}\right)^m
\end{equation}
and
\begin{equation}\label{l}
	T_{e^{il\theta}g}=c_2\left(T_{e^{i\theta}\tilde{\psi}}\right)^l
\end{equation}
for some constants $c_1$ and $c_2$. To avoid the trivial case of having the zero operators, we assume from now on that $c_1$ and $c_2$ are nonzero constants. Finally, using the notion of commutator $[T,S]=TS-ST$ of two operators $T$ and $S$, Equation (3) is then equivalent to
\begin{equation}\label{com}
	c_1\left[\left(T_{e^{i\theta}\widetilde{\phi}}\right)^m,T_{e^{is\theta}r^d}\right](z^k)=c_2\left[T_{e^{ip\theta}r^n},\left(T_{e^{i\theta}\widetilde{\psi}}\right)^l\right](z^k),\textrm{ for all }k\geq 0.
\end{equation}

\section{Preliminaries and tools}
Radial functions in $L^1(\mathbb{D},dA)$ can be seen as functions in $L^1([0,1),rdr)$. For a function $\phi\in L^1([0,1),rdr)$, we define its Mellin transform, denoted $\widehat{\phi}$, by  $$\widehat{\phi}(z)=\int_0^1\phi(r)r^{z-1}\,rdr.$$
It is well known that for those functions  $\phi\in L^1([0,1)rdr)$,  $\widehat{\phi}$ is analytic on the right-half plane $\{z\in\mathbb{C}:\Re z>2\}$ and is continuous and bounded on $\{z\in\mathbb{C}:\Re z\geq2\}$.

Mellin transform (which is related to Laplace transform via the change of variable $r=e^{-t}$) is a precious tool when studying quasihomogeneous Toeplitz operators. In fact, quasihomogeneous Toeplitz operators act on the vectors of the orthogonal basis of $L^2_a(\mathbb{D})$ as shift operators with analytic weight and this weight is involving the Mellin transform of the symbol. We have the following lemma \cite[Lemma 1, p.883]{lry}.
\begin{lemma} Let $\phi\in L^1([0,1), rdr)$ and let $p$ be a non-negative integer. Then for every integer $k\geq 0$, we have 
	$$T_{e^{ip\theta}\phi}(z^k)=2(k+p+1)\widehat{\phi}(2k+p+2)z^{k+p}.$$
\end{lemma}
Another important property of the Mellin transform is that $\widehat{\phi}$ is uniquely determined by its values on any set of integers satisfying the M\"{u}ntz-Sz\'asz (or Blaschke) condition. In our calculations, we
often determine $\phi$ by knowing its Mellin transform $\widehat{\phi}$ on an arithmetic sequences.  We have the following classical theorem \cite[p.102]{Rem}.
\begin{theorem}\label{thm3} 
Suppose that $f$ is a bounded analytic function on the right-half plane $\{z\in\mathbb{C}: \Re z>0\}$ which vanishes at the pairwise distinct points $d_1,d_2,\cdots,$ where
\begin{itemize}
\item[(i)]
$\inf \{\vert d_n\vert\}>0$, and
\item[(ii)]
$\sum\limits_{n\geq 1}^{} \Re(\frac{1}{d_n})=\infty$.
\end{itemize}
Then $f$ vanishes identically on $\{z\in\mathbb{C}:\Re z>0\}$.
\end{theorem}
The next lemma \cite[Lemma 6, p.1468]{l}  is frequently used in the proof of our main result.
\begin{lemma}\label{lem1}
Let $F$ and $G$ be two nonzero bounded analytic functions on the right half-plane $\{z\in\mathbb{C}:\Re (z)>2\}$. If there exists $p\in\mathbb{N}$ such that:
$$
F(z)G(z+p)=F(z+p)G(z),
$$
then $F=c G$, for some constant $c$.
\end{lemma}

The following lemma play a key role in our argument for the proof of the main result. In fact, at a certain stage in the proof, we need to determine when the quotient of four Gamma functions is a rational function. We omit the proof of the lemma, which is a slight modification of \cite[Theorem 3, p.197-198]{cr}. 
\begin{lemma}\label{lem2}
Let $a, b, c, d$ be non-negative integers such that $a+b-c-d=\lambda$ and let $\delta \in\mathbb{N}$. Define the function $H$ to be
\begin{equation*}
H(z)=\frac{\Gamma(\frac{z}{2\delta}+\frac{a}{2\delta})\Gamma(\frac{z}{2\delta}+\frac{b}{2\delta})}{\Gamma(\frac{z}{2\delta}+\frac{c}{2\delta})\Gamma(\frac{z}{2\delta}+\frac{d}{2\delta})}.
\end{equation*}
Then, the function $H$ is a rational function if and only if $2\delta$ divides $\lambda$ and one of the numbers $a-c$ or $a-d$.
 
\end{lemma}
In \cite[Theorem 3, p.197-198]{cr}, the authors assume that $a+b-c-d=-1$ rather than $a+b-c-d=0$ as in our version above. However, the proof remains exactly the same as stated in \cite[p.205]{cr}.  

\section{Main Results}
To ensure clarity and to effectively persuade readers who may not be entirely familiar with the calculations within our proofs, we opt to start by meticulously describing the case when $s$ equals $2p$. This aims to bolster the understanding and conviction that our result holds true. Furthermore, the proof of the general case (Theorem {\ref{thm1}) is essentially based on this key case.

\begin{theorem}\label{thm2}
Let $\phi(r)=r^n$ and $\psi(r)=r^d$ with $p<s$, $n,d\in\mathbb{N}$. Assume, there exist $m, l\in\mathbb{N}$ and nontrivial radial functions $f, g$ such that the hypotheses (H1) and (H2) are satisfied. If $s=2p$, then $m =p$, $l=s$ and
$$
T_{e^{im\theta}f}+T_{e^{is\theta}g}=c(T_{e^{ip\theta}r^n}+ T_{e^{is\theta}r^d}).
$$
for some constant $c$.
\end{theorem}
For the sake of simplicity, we will adopt the following notation in the proof. We shall use $\equiv$ instead of $=$ whenever the quantity on the left side of the equation equals a constant time the quantity on the right side. 
\begin{proof}
Let $p, s, m, l, d, n\in\mathbb{N}^{*}$, and the condition (H1)
\begin{equation}
\Big[\Big(T_{e^{i\theta}f}\Big)^m, T_{e^{is\theta}r^d}\Big](z^k)=\Big[T_{e^{ip\theta}r^n}, \Big(T_{e^{i\theta}g}\Big)^l\Big](z^k),\quad \forall k\in\mathbb{N}.
\end{equation}
such that
\begin{equation}\label{eq5}
T_{e^{ip\theta}r^n}=\Big(T_{e^{i\theta}f}\Big)^p\quad ,\quad T_{e^{is\theta}r^d}= \Big(T_{e^{i\theta}g}\Big)^s.
\end{equation}
then,
\begin{equation}\label{eq1}
\Big[\Big(T_{e^{i\theta}f}\Big)^m, \Big(T_{e^{i\theta}g}\Big)^s\Big](z^k)=\Big[\Big(T_{e^{i\theta}f}\Big)^p, \Big(T_{e^{i\theta}g}\Big)^l\Big](z^k),\quad \forall k\in\mathbb{N}.
\end{equation}
With the condition (H2), namely "$1\leq p<s, 1\leq m<l$ and $m+s=p+l$",
we have
\begin{align*}
&\prod_{j=0}^{s-1}2(k+j+2)\widehat{g}(2k+2j+3)\prod_{j=0}^{m-1}2(k+s+j+2)\widehat{f}(2k+2s+2j+3)\\
&- \prod_{j=0}^{m-1}2(k+j+2)\widehat{f}(2k+2j+3)\prod_{j=0}^{s-1}2(k+m+j+2)\widehat{g}(2k+2m+2j+3)\\
&=\prod_{j=0}^{l-1}2(k+j+2)\widehat{g}(2k+2j+3)\prod_{j=0}^{p-1}2(k+l+j+2)\widehat{f}(2k+2l+2j+3)\\
&- \prod_{j=0}^{p-1}2(k+j+2)\widehat{f}(2k+2j+3)\prod_{j=0}^{l-1}2(k+p+j+2)\widehat{g}(2k+2p+2j+3).
\end{align*}
Thus
\begin{align*}
&\prod_{j=0}^{s-1}2(k+j+2)\widehat{g}(2k+2j+3)\prod_{j=s}^{m+s-1}2(k+j+2)\widehat{f}(2k+2j+3)\\
&- \prod_{j=0}^{m-1}2(k+j+2)\widehat{f}(2k+2j+3)\prod_{j=m}^{m+s-1}2(k+j+2)\widehat{g}(2k+2j+3)\\
&=\prod_{j=0}^{l-1}2(k+j+2)\widehat{g}(2k+2j+3)\prod_{j=l}^{p+l-1}2(k+j+2)\widehat{f}(2k+2j+3)\\
&- \prod_{j=0}^{p-1}2(k+j+2)\widehat{f}(2k+2j+3)\prod_{j=p}^{p+l-1}2(k+j+2)\widehat{g}(2k+2j+3).
\end{align*}
Since $m+s=p+l$, we have that
\begin{align*}
\prod_{j=0}^{s-1}\widehat{g}(2k+2j+3)\prod_{j=s}^{m+s-1}\widehat{f}(2k+2j+3)
&- \prod_{j=0}^{m-1}\widehat{f}(2k+2j+3)\prod_{j=m}^{m+s-1}\widehat{g}(2k+2j+3)\\
=\prod_{j=0}^{l-1}\widehat{g}(2k+2j+3)\prod_{j=l}^{p+l-1}\widehat{f}(2k+2j+3)
&- \prod_{j=0}^{p-1}\widehat{f}(2k+2j+3)\prod_{j=p}^{p+l-1}\widehat{g}(2k+2j+3),
\end{align*}
which yields
\begin{align}\label{eq6}
\prod_{j=0}^{s-1}\widehat{g}(2k+2j+3)\prod_{j=s}^{m+s-1}\widehat{f}(2k+2j+3)
&- \prod_{j=0}^{l-1}\widehat{g}(2k+2j+3)\prod_{j=l}^{p+l-1}\widehat{f}(2k+2j+3)\\
\nonumber =\prod_{j=0}^{m-1}\widehat{f}(2k+2j+3)\prod_{j=m}^{m+s-1}\widehat{g}(2k+2j+3)
&- \prod_{j=0}^{p-1}\widehat{f}(2k+2j+3)\prod_{j=p}^{p+l-1}\widehat{g}(2k+2j+3).
\end{align}
At this stage, we need to study the cases "$l\leq s$ and $m\leq p$" and "$l>s$ and $m>p$" separately.
$${\framebox{{\bf{Case I}}. $l\leq s$ and $m\leq p$.}}$$
In this case, Equation (\ref{eq6}) can be written as
\begin{align*}
\prod_{j=s}^{m+s-1}\widehat{f}(2k+2j+3)\prod_{j=0}^{l-1}\widehat{g}(2k+2j+3)&\Big[\prod_{j=l}^{s-1}\widehat{g}(2k+2j+3)-\prod_{j=l}^{s-1}\widehat{f}(2k+2j+3)\Big]\\
=\prod_{j=0}^{m-1}\widehat{f}(2k+2j+3)\prod_{j=p}^{p+l-1}\widehat{g}(2k+2j+3)&\Big[\prod_{j=m}^{p-1}\widehat{g}(2k+2j+3)-\prod_{j=m}^{p-1}\widehat{f}(2k+2j+3)\Big],
\end{align*}
and so
\begin{equation}\label{eq2}
\dfrac{\prod\limits_{j=s}^{m+s-1}\widehat{f}(2k+2j+3)\prod\limits_{j=0}^{l-1}\widehat{g}(2k+2j+3)}{\prod\limits_{j=0}^{m-1}\widehat{f}(2k+2j+3)\prod\limits_{j=p}^{p+l-1}\widehat{g}(2k+2j+3)}=\dfrac{\prod\limits_{j=m}^{p-1}\widehat{g}(2k+2j+3)-\prod\limits_{j=m}^{p-1}\widehat{f}(2k+2j+3)}{\prod\limits_{j=l}^{s-1}\widehat{g}(2k+2j+3)-\prod\limits_{j=l}^{s-1}\widehat{f}(2k+2j+3)}.
\end{equation}
By (H2) there are $\alpha, \beta\in\mathbb{N}$ such that $m+\beta=l$ and $p+\alpha=s$. Since $l+p=m+s$, we have $m+\beta+p=m+p+\alpha$. Thus we must have $\beta=\alpha$. So let $\alpha\in\mathbb{N}$ be such that $p+\alpha=s, m+\alpha=l$. Then Equation (\ref{eq2}) becomes
\begin{equation}\label{eq3}
\dfrac{\prod\limits_{j=p+\alpha}^{m+s-1}\widehat{f}(2k+2j+3)\prod\limits_{j=0}^{m+\alpha-1}\widehat{g}(2k+2j+3)}{\prod\limits_{j=0}^{m-1}\widehat{f}(2k+2j+3)\prod\limits_{j=p}^{p+m+\alpha-1}\widehat{g}(2k+2j+3)}=\dfrac{\prod\limits_{j=m}^{p-1}\widehat{g}(2k+2j+3)-\prod\limits_{j=m}^{p-1}\widehat{f}(2k+2j+3)}{\prod\limits_{j=m+\alpha}^{p+\alpha-1}\widehat{g}(2k+2j+3)-\prod\limits_{j=m+\alpha}^{p+\alpha-1}\widehat{f}(2k+2j+3)}.
\end{equation}
Since $s=2p$, we have $\alpha=p$. Now, let us define the functions $F, G$ as follows
$$
F(k)=\dfrac{\prod\limits_{j=0}^{m+\alpha-1}\widehat{g}(2k+2j+3)}{\prod\limits_{j=0}^{m-1}\widehat{f}(2k+2j+3)\prod\limits_{j=\alpha}^{m+\alpha-1}\widehat{f}(2k+2j+3)}.
$$
and
$$
G(k)=\prod\limits_{j=m}^{p-1}\widehat{g}(2k+2j+3)-\prod\limits_{j=m}^{p-1}\widehat{f}(2k+2j+3).
$$
It easy to see that Equation (\ref{eq3}) is equivalent to
$$
\frac{F(k)}{F(k+\alpha)}=\frac{G(k)}{G(k+\alpha)},\quad \forall k\in\mathbb{N}.
$$
By Lemma \ref{lem1}, it follows that $G(k)\equiv F(k)$. Thus
\begin{equation}\label{eq4}
\dfrac{\prod\limits_{j=0}^{m+\alpha-1}\widehat{g}(2k+2j+3)}{\prod\limits_{j=0}^{m-1}\widehat{f}(2k+2j+3)\prod\limits_{j=\alpha}^{m+\alpha-1}\widehat{f}(2k+2j+3)}\equiv\prod\limits_{j=m}^{p-1}\widehat{g}(2k+2j+3)-\prod\limits_{j=m}^{p-1}\widehat{f}(2k+2j+3).
\end{equation}
From  equations (\ref{m}) and (\ref{l}), the operator $T_{e^{ip\theta}r^n}$ (resp. $T_{e^{is\theta}r^d}$) commutes with $T_{e^{i\theta}f}$ (resp. $T_{e^{i\theta}g}$). Then, using the techniques in the proof of \cite[Theorem 14, p.1473]{l}, we have for all $k\in\mathbb{N}$
\begin{align*}
\widehat{f}(2k+3)&=B(\frac{2k+4}{2p},1-\frac{1}{p})B(\frac{2k+p+n+2}{2p},\frac{1}{p}) \\
&= \dfrac{\Gamma(\frac{2k+4}{2p})\Gamma(1-\frac{1}{p})}{\Gamma(\frac{2k+2}{2p}+1)}
\dfrac{\Gamma(\frac{2k+p+n+2}{2p})\Gamma(\frac{1}{p})}{\Gamma(\frac{2k+p+n+4}{2p})} \\
&\equiv\dfrac{\Gamma(\frac{2k+4}{2p})\Gamma(\frac{2k+p+n+2}{2p})}{(2k+2)\Gamma(\frac{2k+2}{2p})\Gamma(\frac{2k+p+n+4}{2p})}.
\end{align*}
Similary,
\begin{align*}
\widehat{g}(2k+3)&=B(\frac{2k+4}{2s},1-\frac{1}{s})B(\frac{2k+s+d+2}{2s},\frac{1}{s}) \\
&= \dfrac{\Gamma(\frac{2k+4}{2s})\Gamma(1-\frac{1}{s})}{\Gamma(\frac{2k+2}{2s}+1)}
\dfrac{\Gamma(\frac{2k+s+d+2}{2s})\Gamma(\frac{1}{s})}{\Gamma(\frac{2k+s+d+4}{2s})} \\
&\equiv\dfrac{\Gamma(\frac{2k+4}{2s})\Gamma(\frac{2k+s+d+2}{2s})}{(2k+2)\Gamma(\frac{2k+2}{2s})\Gamma(\frac{2k+s+d+4}{2s})}.
\end{align*}
Thus the terms on the left side of Equation (\ref{eq4}) can be written as
\begin{align*}
\prod\limits_{j=0}^{m+\alpha-1}\widehat{g}(2k+2j+3)&\equiv\prod\limits_{j=0}^{m+\alpha-1}\dfrac{\Gamma(\frac{2k+2j+4}{2s})\Gamma(\frac{2k+2j+s+d+2}{2s})}{(2k+2j+2)\Gamma(\frac{2k+2j+2}{2s})\Gamma(\frac{2k+2j+s+d+4}{2s})} \\
&\equiv \frac{1}{\prod\limits_{j=0}^{m+\alpha-1}(2k+2j+2)}\dfrac{\Gamma(\frac{2k+2m+2\alpha+2}{2s})\Gamma(\frac{2k+s+d+2}{2s})}{\Gamma(\frac{2k+2}{2s})\Gamma(\frac{2k+2m+2\alpha+s+d+2}{2s})} .
\end{align*}
and
\begin{align*}
&\prod\limits_{j=0}^{m-1}\widehat{f}(2k+2j+3)\prod\limits_{j=\alpha}^{m+\alpha-1}\widehat{f}(2k+2j+3) \\
&\equiv \prod\limits_{j=0}^{m-1}\dfrac{\Gamma(\frac{2k+2j+4}{2p})\Gamma(\frac{2k+2j+p+n+2}{2p})}{(2k+2j+2)\Gamma(\frac{2k+2j+2}{2p})\Gamma(\frac{2k+2j+p+n+4}{2p})}\prod\limits_{j=\alpha}^{m+\alpha-1}
\dfrac{\Gamma(\frac{2k+2j+4}{2p})\Gamma(\frac{2k+2j+p+n+2}{2p})}{(2k+2j+2)\Gamma(\frac{2k+2j+2}{2p})\Gamma(\frac{2k+2j+p+n+4}{2p})} \\
&\equiv \prod\limits_{j=0}^{m-1}\prod\limits_{j=\alpha}^{m+\alpha-1}\frac{1}{(2k+2j+2)}
\dfrac{\Gamma(\frac{2k+2m+2}{2p})\Gamma(\frac{2k+p+n+2}{2p})}{\Gamma(\frac{2k+2}{2p})\Gamma(\frac{2k+2m+p+n+2}{2p})}\dfrac{\Gamma(\frac{2k+2m+2\alpha+2}{2p})\Gamma(\frac{2k+2\alpha+p+n+2}{2p})}{\Gamma(\frac{2k+2\alpha+2}{2p})\Gamma(\frac{2k+2m+2\alpha+p+n+2}{2p})} \\
&\equiv \frac{(2k+2m+2)(2k+p+n+2)}{(2k+2)(2k+2m+p+n+2)}\prod\limits_{j=0}^{m-1}\prod\limits_{j=\alpha}^{m+\alpha-1}\frac{1}{(2k+2j+2)}
\Big[\dfrac{\Gamma(\frac{2k+2m+2}{2p})\Gamma(\frac{2k+p+n+2}{2p})}{\Gamma(\frac{2k+2}{2p})\Gamma(\frac{2k+2m+p+n+2}{2p})}\Big]^2,
\end{align*}
so that
\begin{align*}
F(k)&\equiv\frac{(2k+2)(2k+2m+p+n+2)}{(2k+2m+2)(2k+p+n+2)}\prod\limits_{j=0}^{m-1}\prod\limits_{j=\alpha}^{m+\alpha-1}(2k+2j+2)
\Big[\dfrac{\Gamma(\frac{2k+2}{2p})\Gamma(\frac{2k+2m+p+n+2}{2p})}{\Gamma(\frac{2k+2m+2}{2p})\Gamma(\frac{2k+p+n+2}{2p})}\Big]^2 \\
&\times \frac{1}{\prod\limits_{j=0}^{m+\alpha-1}(2k+2j+2)}\dfrac{\Gamma(\frac{2k+2m+2\alpha+2}{2s})\Gamma(\frac{2k+s+d+2}{2s})}{\Gamma(\frac{2k+2}{2s})\Gamma(\frac{2k+2m+2\alpha+s+d+2}{2s})} \\
&\equiv \frac{(2k+2)(2k+2m+p+n+2)}{(2k+2m+2)(2k+p+n+2)}
\Big[\dfrac{\Gamma(\frac{2k+2}{2p})\Gamma(\frac{2k+2m+p+n+2}{2p})}{\Gamma(\frac{2k+2m+2}{2p})\Gamma(\frac{2k+p+n+2}{2p})}\Big]^2 \\
&\times \frac{1}{\prod\limits_{j=m}^{p-1}(2k+2j+2)}\dfrac{\Gamma(\frac{2k+2m+2\alpha+2}{2s})\Gamma(\frac{2k+s+d+2}{2s})}{\Gamma(\frac{2k+2}{2s})\Gamma(\frac{2k+2m+2\alpha+s+d+2}{2s})}.
\end{align*}
On the other hand, the right side of Equation (\ref{eq4}) can be written as
\begin{align*}
&\prod\limits_{j=m}^{p-1}\widehat{g}(2k+2j+3)-\prod\limits_{j=m}^{p-1}\widehat{f}(2k+2j+3) \\
&\equiv\prod\limits_{j=m}^{p-1} \dfrac{\Gamma(\frac{2k+2j+4}{2s})\Gamma(\frac{2k+2j+s+d+2}{2s})}{(2k+2j+2)\Gamma(\frac{2k+2j+2}{2s})\Gamma(\frac{2k+2j+s+d+4}{2s})}-\prod\limits_{j=m}^{p-1}\dfrac{\Gamma(\frac{2k+2j+4}{2p})\Gamma(\frac{2k+2j+p+n+2}{2p})}{(2k+2j+2)\Gamma(\frac{2k+2j+2}{2p})\Gamma(\frac{2k+2j+p+n+4}{2p})} \\
&\equiv\prod\limits_{j=m}^{p-1}\frac{1}{(2k+2j+2)}\Big[ \dfrac{\Gamma(\frac{2k+2p+2}{2s})\Gamma(\frac{2k+2m+s+d+2}{2s})}{\Gamma(\frac{2k+2m+2}{2s})\Gamma(\frac{2k+2p+s+d+2}{2s})}-\dfrac{(2k+2)\Gamma(\frac{2k+2}{2p})\Gamma(\frac{2k+2m+p+n+2}{2p})}{(2k+p+n+2)\Gamma(\frac{2k+2m+2}{2p})\Gamma(\frac{2k+p+n+2}{2p})} \Big].
\end{align*}
By taking $z=2k+2$, Equation (\ref{eq4}) becomes
\begin{align} \label{eq22}
&\frac{z(z+2m+p+n)}{(z+2m)(z+p+n)}
\Big[\dfrac{\Gamma(\frac{z}{2p})\Gamma(\frac{z+2m+p+n}{2p})}{\Gamma(\frac{z+2m}{2p})\Gamma(\frac{z+p+n}{2p})}\Big]^2\dfrac{\Gamma(\frac{z+2m+2\alpha}{2s})\Gamma(\frac{z+s+d}{2s})}{\Gamma(\frac{z}{2s})\Gamma(\frac{z+2m+2\alpha+s+d}{2s})} \\
\nonumber &\equiv \dfrac{\Gamma(\frac{z+2p}{2s})\Gamma(\frac{z+2m+s+d}{2s})}{\Gamma(\frac{z+2m}{2s})\Gamma(\frac{z+2p+s+d}{2s})}-\dfrac{z\Gamma(\frac{z}{2p})\Gamma(\frac{z+2m+p+n}{2p})}{(z+p+n)\Gamma(\frac{z+2m}{2p})\Gamma(\frac{z+p+n}{2p})}.
\end{align}
In (\ref{eq22}), the poles of the function $\Gamma(\frac{z+2m}{2p})\Gamma(\frac{z+p+n}{2p})$ are zeros of the function $\dfrac{\Gamma(\frac{z+2p}{2s})\Gamma(\frac{z+2m+s+d}{2s})}{\Gamma(\frac{z+2m}{2s})\Gamma(\frac{z+2p+s+d}{2s})}$. Thus $\dfrac{\Gamma(\frac{z+2m}{2p})\Gamma(\frac{z+p+n}{2p})}{\Gamma(\frac{z+2m}{2s})\Gamma(\frac{z+2p+s+d}{2s})}$ is a rational function. Also, by taking $z=z+2p$, we can easily show that the function  $\dfrac{\Gamma(\frac{z+2m+2p}{2s})\Gamma(\frac{z+s+d}{2s})}{\Gamma(\frac{z+2m}{2p})\Gamma(\frac{z+p+n}{2p})}$ is rational.
Consequently, the poles of the function $\dfrac{\Gamma(\frac{z+2p}{2s})\Gamma(\frac{z+2m+s+d}{2s})}{\Gamma(\frac{z+2m}{2s})\Gamma(\frac{z+2p+s+d}{2s})}$ are poles of the right side of Equation (\ref{eq22}) and they can come only from the term $\Gamma(\frac{z}{2p})\Gamma(\frac{z+2m+p+n}{2p})$. Thus the function $\dfrac{\Gamma(\frac{z}{2p})\Gamma(\frac{z+2m+p+n}{2p})}{\Gamma(\frac{z+2p}{2s})\Gamma(\frac{z+2m+s+d}{2s})}$ is rational. Again, by taking $z=z+2p$, we can show that the function  $\dfrac{\Gamma(\frac{z}{2p})\Gamma(\frac{z+2m+p+n}{2p})}{\Gamma(\frac{z}{2s})\Gamma(\frac{z+2m+2p+s+d}{2s})}$ is rational. Hence Equation (\ref{eq22}) can be written as
\begin{align} \label{eq23}
&H(z)
\dfrac{\Gamma(\frac{z}{2p})\Gamma(\frac{z+2m+p+n}{2p})}{\Gamma(\frac{z+2m}{2p})\Gamma(\frac{z+p+n}{2p})} \\
\nonumber &\equiv \dfrac{\Gamma(\frac{z+2p}{2s})\Gamma(\frac{z+2m+s+d}{2s})}{\Gamma(\frac{z+2m}{2s})\Gamma(\frac{z+2p+s+d}{2s})}-\dfrac{z\Gamma(\frac{z}{2p})\Gamma(\frac{z+2m+p+n}{2p})}{(z+p+n)\Gamma(\frac{z+2m}{2p})\Gamma(\frac{z+p+n}{2p})},
\end{align}
where $H$ is a rational function. This implies that the functions $\dfrac{\Gamma(\frac{z}{2p})\Gamma(\frac{z+2m+p+n}{2p})}{\Gamma(\frac{z+2m}{2p})\Gamma(\frac{z+p+n}{2p})}$ and $\dfrac{\Gamma(\frac{z+2p}{2s})\Gamma(\frac{z+2m+s+d}{2s})}{\Gamma(\frac{z+2m}{2s})\Gamma(\frac{z+2p+s+d}{2s})}$ are rational functions. Therefore, by Lemma \ref{lem2}, we deduce that $2p$ divides $2m$ or $p+n$ and $2s$ divides $2p-2m$ or $s+d$.
Observe that if $2s$ divides $2p-2m$, then $2p<p-m$ and so $ p<-m<0$ which is a contradiction. That is $2s$ divides $s+d$.
On the other hand if $2p$ divide $2m$, then $p=m$. Hence $s=l$ and this finishes the proof. Finally if $2p$ divides $p+n$, then Theorem \ref{thm4} gives us the desired result.
$${\framebox{{\bf{Case II}}. $l> s$ and $m> p$.}}$$
In this case,  Equation (\ref{eq6}) can be written as
\begin{align*}
\prod_{j=l}^{m+s-1}\widehat{f}(2k+2j+3)\prod_{j=0}^{s-1}\widehat{g}(2k+2j+3)&\Big[\prod_{j=s}^{l-1}\widehat{f}(2k+2j+3)-\prod_{j=s}^{l-1}\widehat{g}(2k+2j+3)\Big]\\
=\prod_{j=0}^{p-1}\widehat{f}(2k+2j+3)\prod_{j=m}^{p+l-1}\widehat{g}(2k+2j+3)&\Big[\prod_{j=p}^{m-1}\widehat{f}(2k+2j+3)-\prod_{j=p}^{m-1}\widehat{g}(2k+2j+3)\Big],
\end{align*}
or
\begin{equation}\label{eq11}
\dfrac{\prod\limits_{j=l}^{m+s-1}\widehat{f}(2k+2j+3)\prod\limits_{j=0}^{s-1}\widehat{g}(2k+2j+3)}{\prod\limits_{j=0}^{p-1}\widehat{f}(2k+2j+3)\prod\limits_{j=m}^{p+l-1}\widehat{g}(2k+2j+3)}=
\dfrac{\prod\limits_{j=p}^{m-1}\widehat{f}(2k+2j+3)-\prod\limits_{j=p}^{m-1}\widehat{g}(2k+2j+3)}{\prod\limits_{j=s}^{l-1}\widehat{f}(2k+2j+3)-\prod\limits_{j=s}^{l-1}\widehat{g}(2k+2j+3)}.
\end{equation}
Multiplying both sides of the above equation by $\dfrac{\prod\limits_{j=s}^{l-1}\widehat{f}(2k+2j+3)\prod\limits_{j=s}^{l-1}\widehat{g}(2k+2j+3)}{\prod\limits_{j=p}^{m-1}\widehat{f}(2k+2j+3)\prod\limits_{j=p}^{m-1}\widehat{g}(2k+2j+3)}$, we obtain
\begin{align}\label{eq7}
& \dfrac{\prod\limits_{j=s}^{m+s-1}\widehat{f}(2k+2j+3)\prod\limits_{j=0}^{l-1}\widehat{g}(2k+2j+3)}{\prod\limits_{j=0}^{m-1}\widehat{f}(2k+2j+3)\prod\limits_{j=p}^{p+l-1}\widehat{g}(2k+2j+3)} \\
&\nonumber =\dfrac{\prod\limits_{j=s}^{l-1}\widehat{f}(2k+2j+3)\prod\limits_{j=s}^{l-1}\widehat{g}(2k+2j+3)}{\prod\limits_{j=p}^{m-1}\widehat{f}(2k+2j+3)\prod\limits_{j=p}^{m-1}\widehat{g}(2k+2j+3)}\Big[\dfrac{\prod\limits_{j=p}^{m-1}\widehat{f}(2k+2j+3)-\prod\limits_{j=p}^{m-1}\widehat{g}(2k+2j+3)}{\prod\limits_{j=s}^{l-1}\widehat{f}(2k+2j+3)-\prod\limits_{j=s}^{l-1}\widehat{g}(2k+2j+3)}\Big].
\end{align}
Similar to the argument in Case I, we show that $p=\alpha$ and we let
$$
G(k)=\dfrac{\prod\limits_{j=p}^{m-1}\widehat{f}(2k+2j+3)-\prod\limits_{j=p}^{m-1}\widehat{g}(2k+2j+3)}{\prod\limits_{j=p}^{m-1}\widehat{f}(2k+2j+3)\prod\limits_{j=p}^{m-1}\widehat{g}(2k+2j+3)}.
$$
Then $G(k)\equiv F(k)$, where $F$ is as in Case I. Thus
\begin{equation}\label{eq8}
\dfrac{\prod\limits_{j=0}^{m+\alpha-1}\widehat{g}(2k+2j+3)}{\prod\limits_{j=0}^{m-1}\widehat{f}(2k+2j+3)\prod\limits_{j=\alpha}^{m+\alpha-1}\widehat{f}(2k+2j+3)}\equiv\dfrac{\prod\limits_{j=p}^{m-1}\widehat{f}(2k+2j+3)-\prod\limits_{j=p}^{m-1}\widehat{g}(2k+2j+3)}{\prod\limits_{j=p}^{m-1}\widehat{f}(2k+2j+3)\prod\limits_{j=p}^{m-1}\widehat{g}(2k+2j+3)}.
\end{equation}
We proceed as in Case I and we write
\begin{align*}
& \prod\limits_{j=p}^{m-1}\widehat{f}(2k+2j+3)\prod\limits_{j=p}^{m-1}\widehat{g}(2k+2j+3) \\
&\equiv \prod\limits_{j=p}^{m-1}\dfrac{\Gamma(\frac{2k+2j+4}{2p})\Gamma(\frac{2k+2j+p+n+2}{2p})}{(2k+2j+2)\Gamma(\frac{2k+2j+2}{2p})\Gamma(\frac{2k+2j+p+n+4}{2p})}\prod\limits_{j=p}^{m-1}\dfrac{\Gamma(\frac{2k+2j+4}{2s})\Gamma(\frac{2k+2j+s+d+2}{2s})}{(2k+2j+2)\Gamma(\frac{2k+2j+2}{2s})\Gamma(\frac{2k+2j+s+d+4}{2s})}  \\
&\equiv \frac{1}{\prod\limits_{j=p}^{m-1}(2k+2j+2)^2}\dfrac{\Gamma(\frac{2k+2m+2}{2p})\Gamma(\frac{2k+2p+p+n+2}{2p})}{\Gamma(\frac{2k+2p+2}{2p})\Gamma(\frac{2k+2m+p+n+2}{2p})}\dfrac{\Gamma(\frac{2k+2m+2}{2s})\Gamma(\frac{2k+2p+s+d+2}{2s})}{\Gamma(\frac{2k+2p+2}{2s})\Gamma(\frac{2k+2m+s+d+2}{2s})}  \\
&\equiv \frac{1}{\prod\limits_{j=p}^{m-1}(2k+2j+2)^2}\dfrac{(2k+p+n+2)\Gamma(\frac{2k+2m+2}{2p})\Gamma(\frac{2k+p+n+2}{2p})}{(2k+2)\Gamma(\frac{2k+2}{2p})\Gamma(\frac{2k+2m+p+n+2}{2p})}\dfrac{\Gamma(\frac{2k+2m+2}{2s})\Gamma(\frac{2k+2p+s+d+2}{2s})}{\Gamma(\frac{2k+2p+2}{2s})\Gamma(\frac{2k+2m+s+d+2}{2s})}.
\end{align*}
and
\begin{align*}
&\prod\limits_{j=p}^{m-1}\widehat{f}(2k+2j+3)-\prod\limits_{j=p}^{m-1}\widehat{g}(2k+2j+3) \\
&\equiv\prod\limits_{j=p}^{m-1}\dfrac{\Gamma(\frac{2k+2j+4}{2p})\Gamma(\frac{2k+2j+p+n+2}{2p})}{(2k+2j+2)\Gamma(\frac{2k+2j+2}{2p})\Gamma(\frac{2k+2j+p+n+4}{2p})}-\prod\limits_{j=p}^{m-1}
\dfrac{\Gamma(\frac{2k+2j+4}{2s})\Gamma(\frac{2k+2j+s+d+2}{2s})}{(2k+2j+2)\Gamma(\frac{2k+2j+2}{2s})\Gamma(\frac{2k+2j+s+d+4}{2s})}\\
&\equiv\prod\limits_{j=p}^{m-1}\frac{1}{(2k+2j+2)}\Big[\dfrac{\Gamma(\frac{2k+2m+2}{2p})\Gamma(\frac{2k+2p+p+n+2}{2p})}{\Gamma(\frac{2k+2p+2}{2p})\Gamma(\frac{2k+2m+p+n+2}{2p})}-
\dfrac{\Gamma(\frac{2k+2m+2}{2s})\Gamma(\frac{2k+2p+s+d+2}{2s})}{\Gamma(\frac{2k+2p+2}{2s})\Gamma(\frac{2k+2m+s+d+2}{2s})} \Big] \\
&\equiv\prod\limits_{j=p}^{m-1}\frac{1}{(2k+2j+2)}\Big[\dfrac{(2k+p+n+2)\Gamma(\frac{2k+2m+2}{2p})\Gamma(\frac{2k+p+n+2}{2p})}{(2k+2)\Gamma(\frac{2k+2}{2p})\Gamma(\frac{2k+2m+p+n+2}{2p})}-
\dfrac{\Gamma(\frac{2k+2m+2}{2s})\Gamma(\frac{2k+2p+s+d+2}{2s})}{\Gamma(\frac{2k+2p+2}{2s})\Gamma(\frac{2k+2m+s+d+2}{2s})} \Big].
\end{align*}
By taking $z=2k+2$, Equation (\ref{eq8}) becomes
\begin{align}\label{eq9}
& \frac{(z+2m+p+n)}{(z+2m)}
\dfrac{\Gamma(\frac{z}{2p})\Gamma(\frac{z+2m+p+n}{2p})}{\Gamma(\frac{z+2m}{2p})\Gamma(\frac{z+p+n}{2p})} \dfrac{\Gamma(\frac{z+2m+2\alpha}{2s})\Gamma(\frac{z+s+d}{2s})}{\Gamma(\frac{z}{2s})\Gamma(\frac{z+2m+2\alpha+s+d}{2s})}
\dfrac{\Gamma(\frac{z+2m}{2s})\Gamma(\frac{z+2p+s+d}{2s})}{\Gamma(\frac{z+2p}{2s})\Gamma(\frac{z+2m+s+d}{2s})} \\
&\nonumber\equiv \dfrac{(z+p+n)\Gamma(\frac{z+2m}{2p})\Gamma(\frac{z+p+n}{2p})}{z\Gamma(\frac{z}{2p})\Gamma(\frac{z+2m+p+n}{2p})}-\dfrac{\Gamma(\frac{z+2m}{2s})\Gamma(\frac{z+2p+s+d}{2s})}{\Gamma(\frac{z+2p}{2s})\Gamma(\frac{z+2m+s+d}{2s})} .
\end{align}
Clearly, the poles of the function $\Gamma(\frac{z+2p}{2s})\Gamma(\frac{z+2m+s+d}{2s})$ are zeros of the right side of Equation (\ref{eq9}). Thus $\frac{\Gamma(\frac{z}{2p})\Gamma(\frac{z+2m+p+n}{2p})}{\Gamma(\frac{z+2p}{2s})\Gamma(\frac{z+2m+s+d}{2s})}$ must be a rational function. Consequently, the poles of the function $\Gamma(\frac{z+2m}{2p})\Gamma(\frac{z+p+n}{2p})$ in the right side of Equation (\ref{eq9}) can come only from the term $\Gamma(\frac{z+2m+2\alpha}{2s})\Gamma(\frac{z+s+d}{2s})$, and hence the function $\frac{\Gamma(\frac{z+2m+2\alpha}{2s})\Gamma(\frac{z+s+d}{2s})}{\Gamma(\frac{z+2m}{2p})\Gamma(\frac{z+p+n}{2p})}$ is rational.
Now, we can write  Equation (\ref{eq9}) as follows
\begin{align}\label{eq10}
& H(z)
\dfrac{\Gamma(\frac{z+2m}{2s})\Gamma(\frac{z+2p+s+d}{2s})}{\Gamma(\frac{z}{2s})\Gamma(\frac{z+2m+2\alpha+s+d}{2s})} \\
&\nonumber= \dfrac{(z+p+n)\Gamma(\frac{z+2m}{2p})\Gamma(\frac{z+p+n}{2p})}{z\Gamma(\frac{z}{2p})\Gamma(\frac{z+2m+p+n}{2p})}-\dfrac{\Gamma(\frac{z+2m}{2s})\Gamma(\frac{z+2p+s+d}{2s})}{\Gamma(\frac{z+2p}{2s})\Gamma(\frac{z+2m+s+d}{2s})},
\end{align}
with $H$ being a rational function. In addition, the function $\frac{\Gamma(\frac{z+2m}{2p})\Gamma(\frac{z+p+n}{2p})}{\Gamma(\frac{z}{2p})\Gamma(\frac{z+2m+p+n}{2p})}$ must be rational too. Hence we deduce from Lemma \ref{lem2} that $2p$ divides $2m$ or $p+n$. 
Taking $z=z+2p$ in the function $\frac{\Gamma(\frac{z+2m+2\alpha}{2s})\Gamma(\frac{z+s+d}{2s})}{\Gamma(\frac{z+2m}{2p})\Gamma(\frac{z+p+n}{2p})}$ and multiplying the result by  $\frac{\Gamma(\frac{z}{2p})\Gamma(\frac{z+2m+p+n}{2p})}{\Gamma(\frac{z+2p}{2s})\Gamma(\frac{z+2m+s+d}{2s})}$, we obtain that $\frac{\Gamma(\frac{z+2m}{2s})\Gamma(\frac{z+2p+s+d}{2s})}{\Gamma(\frac{z+2p}{2s})\Gamma(\frac{z+2m+s+d}{2s})}$ is a rational function.
From Equation (\ref{eq9}), we find that the functions $\frac{\Gamma(\frac{z+2m}{2s})\Gamma(\frac{z+2p+s+d}{2s})}{\Gamma(\frac{z}{2s})\Gamma(\frac{z+2m+2\alpha+s+d}{2s})}$ and $\frac{\Gamma(\frac{z+2m}{2s})\Gamma(\frac{z+2p+s+d}{2s})}{\Gamma(\frac{z+2p}{2s})\Gamma(\frac{z+2m+s+d}{2s})}$ should be rational as well.
More precisely, by Lemma \ref{lem2} we have the following possibilities: $2s$ divides $2m$ or $2p+s+d$, $2s$ divides $2m-2p$ or $-s-d$, and $2p$ divides $2m$ or $p+n$.
If $\frac{m}{p}\notin\mathbb{N}$, then it's straightforward to see that $2s$ divides $s+d$ and $2p$ divides $p+n$, and so the conclusion follows from Theorem \ref{thm4}.
We now focus on the case when $\frac{m}{p}\in\mathbb{N}$, and we discuss the following two situations.
{\bf{Situation 1}}: $\mathbf{m=(2N+1)p}$, where $\mathbf{N\in\mathbb{N}}$.
Equation (\ref{eq9}) becomes
\begin{align}\label{eq13}
& \frac{(z+2m+p+n)}{(z+2m)}
\prod_{i=0}^{2N}\frac{(z+p+n+2pi)}{(z+2pi)}\prod_{i=0}^{N}\frac{(z+2si)}{(z+s+d+2si)}
\prod_{i=0}^{N-1}\frac{(z+2p+2si)}{(z+2p+s+d+2si)} \\
&\nonumber \equiv \prod_{i=1}^{2N}\frac{(z+2pi)}{(z+p+n+2pi)}-\prod_{i=0}^{N-1}\frac{(z+2p+2si)}{(z+2p+s+d+2si)}.
\end{align}
Observe that the points $-p-n-2pi$ for all $1\leq i\leq 2N$ are poles of the right side. But in the same time those points are zeros of the left side. So some of them must be canceled, which means that $2p$ divides $p+n$. Therefore there exists $1\leq t_1\leq 2N-1$ such that $2pt_1=p+n$. In this case, Equation (\ref{eq13}) becomes
\begin{align}\label{eq15}
& \frac{(z+2m+p+n)}{(z+2m)}
\frac{\prod\limits_{i=2N+1}^{2N+t_1}(z+2pi)}{\prod\limits_{i=0}^{t_1-1}(z+2pi)}\prod_{i=0}^{N}\frac{(z+2si)}{(z+s+d+2si)}
\prod_{i=0}^{N-1}\frac{(z+2p+2si)}{(z+2p+s+d+2si)} \\
&\nonumber \equiv \frac{\prod\limits_{i=1}^{t_1}(z+2pi)}{\prod\limits_{i=2N+1}^{t_1+2N}(z+2pi)}-\prod_{i=0}^{N-1}\frac{(z+2p+2si)}{(z+2p+s+d+2si)}.
\end{align}
We shall argue by contradiction. Let us assume that $2s$ cannot divide $s+d$. Then, the second product on the right side cannot be reduced or simplified.
In this situation, we have two subcases:
\begin{itemize}
\item[(i)]
If $2st_2=2p+s+d$ with $1\leq t_2\leq N$, then Equation (\ref{eq15}) becomes
\begin{align*}
& \frac{(z+2m+p+n)}{(z+2m)}
\frac{\prod\limits_{i=2N+1}^{t_1+2N}(z+2pi)}{\prod\limits_{i=0}^{t_1-1}(z+2pi)}\frac{\prod\limits_{i=0}^{N}(z+2si)}{\prod\limits_{i=t_2-1}^{t_2+N-1}(z+2si+2p)}
\frac{\prod\limits_{i=0}^{N-1}(z+2p+2si)}{\prod\limits_{i=t_2}^{t_2+N-1}(z+2si)} \\
&\equiv \frac{\prod\limits_{i=1}^{t_1}(z+2pi)}{\prod\limits_{i=2N+1}^{t_1+2N}(z+2pi)}-\frac{\prod\limits_{i=0}^{N-1}(z+2p+2si)}{\prod\limits_{i=t_2}^{t_2+N-1}(z+2si)},
\end{align*}
or
\begin{align*}
& \frac{(z+2m+p+n)}{(z+2m)}
\frac{\prod\limits_{i=2N+1}^{t_1+2N}(z+2pi)}{\prod\limits_{i=0}^{t_1-1}(z+2pi)}\frac{\prod\limits_{i=0}^{t_2-1}(z+2si)}{\prod\limits_{i=N}^{t_2+N-1}(z+2si+2p)}
\frac{\prod\limits_{i=0}^{t_2-2}(z+2p+2si)}{\prod\limits_{i=N+1}^{t_2+N-1}(z+2si)} \\
&\equiv \frac{\prod\limits_{i=1}^{t_1}(z+2pi)}{\prod\limits_{i=2N+1}^{t_1+2N}(z+2pi)}-\frac{\prod\limits_{i=0}^{N-1}(z+2p+2si)}{\prod\limits_{i=t_2}^{t_2+N-1}(z+2si)} .
\end{align*}
Now the point $-2sN$ is pole of the right side, but on the left side it can come only from $-2pi$ where $0\leq i\leq t_1-1$. However this cannot be true because $2pi\leq 2p(2N-2)=2sN-2s$. Hence $2s$ must divide $s+d$ which contradicts our assumption.
\item[(ii)]
If $2s$ cannot divided $2p+s+d$, then, from Equation (\ref{eq15}), the poles $-s-d-2sj$ with $0\leq j\leq N$ in the left side of Equation (\ref{eq15}) can appear on the right side only from $-2pi$ with $2N+1\leq i\leq t_1+2N$.
Observe that it's easy to get the contradiction here if $t_1= 1$. Now, we set $i=2N+2$ and we obtain $2p(2N+2)=4pN+4p=2sN+2s=s+d+2sj$. Hence $2s$ divides $s+d$, which is a contradiction.
\end{itemize}
{\bf{Situation 2}}: $\mathbf{m=2Np}$, where $\mathbf{N\in\mathbb{N}}$.
In this case, Equation (\ref{eq9}) becomes 
\begin{align*}
& \frac{(z+2m+p+n)}{(z+2m)}\prod_{i=0}^{2N-1}\frac{(z+p+n+2pi)}{(z+2pi)}\prod_{i=0}^{N-1}\frac{(z+2si)}{(z+2p+s+d+2si)}\prod_{i=0}^{N-1}\frac{(z+2p+2si)}{(z+s+d+2si)}   \\
&\equiv \prod_{i=1}^{2N-1}\frac{(z+2pi)}{(z+p+n+2pi)}-\dfrac{\Gamma(\frac{z+2m}{2s})\Gamma(\frac{z+2p+s+d}{2s})}{\Gamma(\frac{z+2p}{2s})\Gamma(\frac{z+2m+s+d}{2s})} .
\end{align*} 
We notice here that $2s$ cannot divide $2m-2p$, and so it divides $s+d$. We can then write that $2st_3=s+d$ for some $t_3\in\mathbb{N}$. Consequently, the previous equation can be written as
\begin{align}\label{eq16}
& \frac{(z+2m+p+n)}{(z+2m)}\prod_{i=0}^{2N-1}\frac{(z+p+n+2pi)}{(z+2pi)}\frac{\prod\limits_{i=0}^{N-1}(z+2si)}{\prod\limits_{i=t_3}^{N+t_3-1}(z+2p+2si)}\frac{\prod\limits_{i=0}^{N-1}(z+2p+2si)}{\prod\limits_{i=t_3}^{N+t_3-1}(z+2si)}   \\
&\nonumber \equiv \prod_{i=1}^{2N-1}\frac{(z+2pi)}{(z+p+n+2pi)}-\frac{\prod\limits_{i=0}^{t_3-1}(z+2p+2si)}{\prod\limits_{i=N}^{N+t_3-1}(z+2si)} .
\end{align}
We observe that the points $-p-n-2pi$ are poles of the right side of Equation (\ref{eq16}) for all $1\leq i\leq 2N-1$. But this cannot hold because all the poles of the left side are multiple of $2p$. Hence  we necessary have that $2p$ divides $p+n$ and the proof is complete.
\end{proof}
For the general case, we shall drop the condition we had for Theorem 3, namely $s=2p$. As mentioned in the beginning of this section, the proof of the next result is closely related to the calculations done in the proof of Theorem 3.
\begin{theorem}\label{thm1}
Let $\phi(r)=r^n$ and $\psi(r)=r^d$ with $p<s$, $n,d\in\mathbb{N}$. Assume, there exist $m, l\in\mathbb{N}$ and nontrivial radial functions $f, g$ such that the hypotheses (H1) and (H2) are satisfied. Then $m =p$, $l=s$ and
$$
T_{e^{im\theta}f}+T_{e^{is\theta}g}=c(T_{e^{ip\theta}r^n}+ T_{e^{is\theta}r^d}),
$$
for some constant $c$.
\end{theorem}
\begin{proof}
As in the proof of Theorem (\ref{thm2}), we treat the cases "$l\leq s$ and $m\leq p$" and $l>s$ and $m> p$" separately.
$${\framebox{{\bf{Case I}}. $l\leq s$ and $m\leq p$.}}$$
 Similarly as in Case I of the proof of Theorem 3,  Equation (\ref{eq3}) becomes
\begin{equation}\label{eq33}
\dfrac{\prod\limits_{j=p+\alpha}^{m+s-1}\widehat{f}(2k+2j+3)\prod\limits_{j=0}^{m+\alpha-1}\widehat{g}(2k+2j+3)}{\prod\limits_{j=0}^{m-1}\widehat{f}(2k+2j+3)\prod\limits_{j=p}^{p+m+\alpha-1}\widehat{g}(2k+2j+3)}=\dfrac{\prod\limits_{j=m}^{p-1}\widehat{g}(2k+2j+3)-\prod\limits_{j=m}^{p-1}\widehat{f}(2k+2j+3)}{\prod\limits_{j=m+\alpha}^{p+\alpha-1}\widehat{g}(2k+2j+3)-\prod\limits_{j=m+\alpha}^{p+\alpha-1}\widehat{f}(2k+2j+3)},
\end{equation}
or
\begin{align*}
& \dfrac{\prod\limits_{j=p+\alpha}^{m+s-1}\widehat{f}(2k+2j+3)\prod\limits_{j=0}^{m+\alpha-1}\widehat{g}(2k+2j+3)}{\prod\limits_{j=0}^{m-1}\widehat{f}(2k+2j+3)\prod\limits_{j=p}^{p+m+\alpha-1}\widehat{g}(2k+2j+3)}  \\
&\equiv \prod\limits_{j=p+\alpha}^{m+s-1}\dfrac{\Gamma(\frac{2k+2j+4}{2p})\Gamma(\frac{2k+2j+p+n+2}{2p})}{(2k+2j+2)\Gamma(\frac{2k+2j+2}{2p})\Gamma(\frac{2k+2j+p+n+4}{2p})}
\prod\limits_{j=0}^{m-1} \dfrac{(2k+2j+2)\Gamma(\frac{2k+2j+2}{2p})\Gamma(\frac{2k+2j+p+n+4}{2p})}{\Gamma(\frac{2k+2j+4}{2p})\Gamma(\frac{2k+2j+p+n+2}{2p})} \\
&\times
\prod\limits_{j=0}^{m+\alpha-1} \dfrac{\Gamma(\frac{2k+2j+4}{2s})\Gamma(\frac{2k+2j+s+d+2}{2s})}{(2k+2j+2)\Gamma(\frac{2k+2j+2}{2s})\Gamma(\frac{2k+2j+s+d+4}{2s})}
\prod\limits_{j=p}^{p+m+\alpha-1} \dfrac{(2k+2j+2)\Gamma(\frac{2k+2j+2}{2s})\Gamma(\frac{2k+2j+s+d+4}{2s})}{\Gamma(\frac{2k+2j+4}{2s})\Gamma(\frac{2k+2j+s+d+2}{2s})}  \\
&\equiv\frac{\prod\limits_{j=0}^{m-1}(2k+2j+2)}{\prod\limits_{j=p+\alpha}^{m+s-1}(2k+2j+2)} \dfrac{\Gamma(\frac{2k+2m+2s+2}{2p})\Gamma(\frac{2k+2p+2\alpha+p+n+2}{2p})}{\Gamma(\frac{2k+2p+2\alpha+2}{2p})\Gamma(\frac{2k+2m+2s+p+n+2}{2p})}
\dfrac{\Gamma(\frac{2k+2}{2p})\Gamma(\frac{2k+2m+p+n+2}{2p})}{\Gamma(\frac{2k+2m+2}{2p})\Gamma(\frac{2k+p+n+2}{2p})} \\
& \times
\frac{\prod\limits_{j=p}^{s+m-1}(2k+2j+2)}{\prod\limits_{j=0}^{m+\alpha-1}(2k+2j+2)}\dfrac{\Gamma(\frac{2k+2m+2\alpha+2}{2s})\Gamma(\frac{2k+s+d+2}{2s})}{\Gamma(\frac{2k+2}{2s})\Gamma(\frac{2k+2m+2\alpha+s+d+2}{2s})}
\dfrac{\Gamma(\frac{2k+2p+2}{2s})\Gamma(\frac{2k+2m+2s+s+d+2}{2s})}{\Gamma(\frac{2k+2m+2s+2}{2s})\Gamma(\frac{2k+2p+s+d+2}{2s})} \\
&\equiv \dfrac{\Gamma(\frac{2k+2m+2s+2}{2p})\Gamma(\frac{2k+2p+2\alpha+p+n+2}{2p})}{\Gamma(\frac{2k+2p+2\alpha+2}{2p})\Gamma(\frac{2k+2m+2s+p+n+2}{2p})}
\dfrac{\Gamma(\frac{2k+2}{2p})\Gamma(\frac{2k+2m+p+n+2}{2p})}{\Gamma(\frac{2k+2m+2}{2p})\Gamma(\frac{2k+p+n+2}{2p})} \\
& \times
\frac{\prod\limits_{j=p}^{p+\alpha-1}(2k+2j+2)}{\prod\limits_{j=m}^{m+\alpha-1}(2k+2j+2)}\dfrac{\Gamma(\frac{2k+2m+2\alpha+2}{2s})\Gamma(\frac{2k+s+d+2}{2s})}{\Gamma(\frac{2k+2}{2s})\Gamma(\frac{2k+2m+2\alpha+s+d+2}{2s})}
\dfrac{\Gamma(\frac{2k+2p+2}{2s})\Gamma(\frac{2k+2m+2s+s+d+2}{2s})}{\Gamma(\frac{2k+2m+2s+2}{2s})\Gamma(\frac{2k+2p+s+d+2}{2s})} .
\end{align*}
By taking $z=2k+2$, the left side "LS" becomes
\begin{align*}
LS &\equiv \dfrac{(z+2l)(z+2\alpha+p+n)\Gamma(\frac{z+2l}{2p})\Gamma(\frac{z+2\alpha+p+n}{2p})}{(z+2\alpha)(z+2l+p+n)\Gamma(\frac{z+2\alpha}{2p})\Gamma(\frac{z+2l+p+n}{2p})}
\dfrac{\Gamma(\frac{z}{2p})\Gamma(\frac{z+2m+p+n}{2p})}{\Gamma(\frac{z+2m}{2p})\Gamma(\frac{z+p+n}{2p})} \\
& \times
\frac{\prod\limits_{j=p}^{s-1}(z+2j)}{\prod\limits_{j=m}^{l-1}(z+2j)}\dfrac{(z+2m+s+d)\Gamma(\frac{z+2m+2\alpha}{2s})\Gamma(\frac{z+s+d}{2s})}{(z+2m)\Gamma(\frac{z}{2s})\Gamma(\frac{z+2m+2\alpha+s+d}{2s})}
\dfrac{\Gamma(\frac{z+2p}{2s})\Gamma(\frac{z+2m+s+d}{2s})}{\Gamma(\frac{z+2m}{2s})\Gamma(\frac{z+2p+s+d}{2s})}   \\
&\equiv \dfrac{(z+2l)(z+2\alpha+p+n)(z+2m+s+d)(z)\prod\limits_{j=p}^{s-1}(z+2j)}{(z+2\alpha)(z+2l+p+n)(z+2m)(z+s+d)\prod\limits_{j=m}^{l-1}(z+2j)}\frac{G(z+2\alpha)}{G(z)},
\end{align*}
where 
\begin{equation*}
G(z)=\frac{\Gamma(\frac{z+2m}{2p})\Gamma(\frac{z+p+n}{2p})}{\Gamma(\frac{z}{2p})\Gamma(\frac{z+2m+p+n}{2p})}\frac{\Gamma(\frac{z+2m}{2s})\Gamma(\frac{z+2p+s+d}{2s})}{\Gamma(\frac{z+2p}{2s})\Gamma(\frac{z+2m+s+d}{2s})}.
\end{equation*}
Now, we have
\begin{align*}
& \prod\limits_{j=m}^{p-1}\widehat{g}(2k+2j+3)-\prod\limits_{j=m}^{p-1}\widehat{f}(2k+2j+3) \\
&\equiv \prod\limits_{j=m}^{p-1} \dfrac{\Gamma(\frac{2k+2j+4}{2s})\Gamma(\frac{2k+2j+s+d+2}{2s})}{(2k+2j+2)\Gamma(\frac{2k+2j+2}{2s})\Gamma(\frac{2k+2j+s+d+4}{2s})} -\prod\limits_{j=m}^{p-1}
\dfrac{\Gamma(\frac{2k+2j+4}{2p})\Gamma(\frac{2k+2j+p+n+2}{2p})}{(2k+2j+2)\Gamma(\frac{2k+2j+2}{2p})\Gamma(\frac{2k+2j+p+n+4}{2p})}  \\
&\equiv\frac{1}{\prod\limits_{j=m}^{p-1}(z+2j)}\Big[\dfrac{\Gamma(\frac{z+2p}{2s})\Gamma(\frac{z+2m+s+d}{2s})}{\Gamma(\frac{z+2m}{2s})\Gamma(\frac{z+2p+s+d}{2s})}- \dfrac{(z)\Gamma(\frac{z}{2p})\Gamma(\frac{z+2m+p+n}{2p})}{(z+p+n)\Gamma(\frac{z+2m}{2p})\Gamma(\frac{z+p+n}{2p})}\Big] .
\end{align*}
So Equation (\ref{eq33}) becomes
\begin{equation}\label{eq21}
H(z)F(z+2\alpha)\equiv F(z),
\end{equation}
with
\begin{align*}
F(z) &= G(z)\dfrac{\Gamma(\frac{z+2p}{2s})\Gamma(\frac{z+2m+s+d}{2s})}{\Gamma(\frac{z+2m}{2s})\Gamma(\frac{z+2p+s+d}{2s})}- G(z)\dfrac{(z)\Gamma(\frac{z}{2p})\Gamma(\frac{z+2m+p+n}{2p})}{(z+p+n)\Gamma(\frac{z+2m}{2p})\Gamma(\frac{z+p+n}{2p})} \\
&= \frac{z+2m}{z}\Big[\frac{\Gamma(\frac{z+2m}{2p})\Gamma(\frac{z+p+n}{2p})}{\Gamma(\frac{z}{2p})\Gamma(\frac{z+2m+p+n}{2p})}-\frac{(z)\Gamma(\frac{z+2m}{2s})\Gamma(\frac{z+2p+s+d}{2s})}{(z+p+n)\Gamma(\frac{z+2p}{2s})\Gamma(\frac{z+2m+s+d}{2s})}\Big]
\end{align*}
and 
\begin{align*}
H(z) &= \dfrac{(z+2l)(z+2\alpha+p+n)(z+2m+s+d)(z)\prod\limits_{j=p}^{s-1}(z+2j)}{(z+2\alpha)(z+2l+p+n)(z+2m)(z+s+d)\prod\limits_{j=m}^{l-1}(z+2j)}\times \frac{\prod\limits_{j=m}^{p-1}(z+2j)}{\prod\limits_{j=m+\alpha}^{p+\alpha-1}(z+2j)} \\
&= \dfrac{(z+2\alpha+p+n)(z+2m+s+d)}{(z+2l+p+n)(z+s+d)} .
\end{align*}
Our aim now is to show that the function $F$ has only finitely many poles. In fact $F$ has poles at the points $-2sA-2m, -2sB-2p-s-d, -2pC-2m, -2pD-p-n$ for all large $A, B, C, D\in\mathbb{N}$ and These poles come from the terms $\Gamma(\frac{z+2m}{2s})\Gamma(\frac{z+2p+s+d}{2s})$ and $\Gamma(\frac{z+2m}{2p})\Gamma(\frac{z+p+n}{2p})$. In addition, Equation (\ref{eq21}) implies that those points is also poles of the function $F(z+2\alpha)$. We then have the following two situations.
\begin{itemize}
\item[\textbf{Situation 1:}]
The poles of $\Gamma(\frac{z+2m}{2s})\Gamma(\frac{z+2p+s+d}{2s})$ come from $\Gamma(\frac{z+2\alpha+2m}{2s})\Gamma(\frac{z+s+d}{2s})$ and the poles of $\Gamma(\frac{z+2m}{2p})\Gamma(\frac{z+p+n}{2p})$ come from $\Gamma(\frac{z+2\alpha+2m}{2p})\Gamma(\frac{z+2\alpha+p+n}{2p})$, which means that the functions $\frac{\Gamma(\frac{z+2m}{2s})\Gamma(\frac{z+2p+s+d}{2s})}{\Gamma(\frac{z+2\alpha+2m}{2s})\Gamma(\frac{z+s+d}{2s})}$ and $\frac{\Gamma(\frac{z+2m}{2p})\Gamma(\frac{z+p+n}{2p})}{\Gamma(\frac{z+2\alpha+2m}{2p})\Gamma(\frac{z+2\alpha+p+n}{2p})}$ are rational. Hence Lemma \ref{lem2} implies that $s$ divides $p-\alpha$ and $p$ divides $2\alpha$. Therefore $p=\alpha$ and Theorem \ref{thm2} finishes the proof. 
\item[\textbf{Situation 2:}]
The poles of $\Gamma(\frac{z+2m}{2s})\Gamma(\frac{z+2p+s+d}{2s})$ come from $\Gamma(\frac{z+2\alpha+2m}{2p})\Gamma(\frac{z+2\alpha+p+n}{2p})$ and the poles of $\Gamma(\frac{z+2m}{2p})\Gamma(\frac{z+p+n}{2p})$ come from $\Gamma(\frac{z+2\alpha+2m}{2s})\Gamma(\frac{z+s+d}{2s})$, which means that the functions $\frac{\Gamma(\frac{z+2m}{2s})\Gamma(\frac{z+2p+s+d}{2s})}{\Gamma(\frac{z+2\alpha+2m}{2p})\Gamma(\frac{z+2\alpha+p+n}{2p})}$ and $\frac{\Gamma(\frac{z+2m}{2p})\Gamma(\frac{z+p+n}{2p})}{\Gamma(\frac{z+2\alpha+2m}{2s})\Gamma(\frac{z+s+d}{2s})}$ are rational.
Moreover, the function $\frac{\Gamma(\frac{z+2m}{2p})\Gamma(\frac{z+p+n}{2p})}{\Gamma(\frac{z+2m}{2s})\Gamma(\frac{z+2p+s+d}{2s})}$ is also rational. So when  multiplying it by $\frac{\Gamma(\frac{z+2\alpha+2m}{2s})\Gamma(\frac{z+s+d}{2s})}{\Gamma(\frac{z+2m}{2p})\Gamma(\frac{z+p+n}{2p})}$ we obtain that the function $\frac{\Gamma(\frac{z+2\alpha+2m}{2s})\Gamma(\frac{z+s+d}{2s})}{\Gamma(\frac{z+2m}{2s})\Gamma(\frac{z+2p+s+d}{2s})}$ is rational, which implies by Lemma \ref{lem2} that $s$ divides $\alpha-p$. Now if $\alpha>p$, then $s=p+\alpha\leq \alpha -p$ which is not possible. Finally if $\alpha< p$, then $s=p+\alpha\leq p-\alpha$ which is also not possible. Therefore $p=\alpha$  and again Theorem \ref{thm2} finishes the proof.
\end{itemize}
$${\framebox{{\bf{Case II}}. $l> s$ and $m> p$.}}$$
In this case, we can write Equation (\ref{eq6}) as
\begin{align*}
\prod_{j=l}^{m+s-1}\widehat{f}(2k+2j+3)\prod_{j=0}^{s-1}\widehat{g}(2k+2j+3)&\Big[\prod_{j=s}^{l-1}\widehat{f}(2k+2j+3)-\prod_{j=s}^{l-1}\widehat{g}(2k+2j+3)\Big]\\
=\prod_{j=0}^{p-1}\widehat{f}(2k+2j+3)\prod_{j=m}^{p+l-1}\widehat{g}(2k+2j+3)&\Big[\prod_{j=p}^{m-1}\widehat{f}(2k+2j+3)-\prod_{j=p}^{m-1}\widehat{g}(2k+2j+3)\Big]
\end{align*}
or
\begin{equation}\label{eq20}
\dfrac{\prod\limits_{j=l}^{m+s-1}\widehat{f}(2k+2j+3)\prod\limits_{j=0}^{s-1}\widehat{g}(2k+2j+3)}{\prod\limits_{j=0}^{p-1}\widehat{f}(2k+2j+3)\prod\limits_{j=m}^{p+l-1}\widehat{g}(2k+2j+3)}=
\dfrac{\prod\limits_{j=p}^{m-1}\widehat{f}(2k+2j+3)-\prod\limits_{j=p}^{m-1}\widehat{g}(2k+2j+3)}{\prod\limits_{j=s}^{l-1}\widehat{f}(2k+2j+3)-\prod\limits_{j=s}^{l-1}\widehat{g}(2k+2j+3)}.
\end{equation}
Since
\begin{align*}
 & \dfrac{\prod\limits_{j=l}^{m+s-1}\widehat{f}(2k+2j+3)}{\prod\limits_{j=m}^{p+l-1}\widehat{g}(2k+2j+3)} \\
&\equiv\prod\limits_{j=l}^{m+s-1}\dfrac{\Gamma(\frac{2k+2j+4}{2p})\Gamma(\frac{2k+2j+p+n+2}{2p})}{(2k+2j+2)\Gamma(\frac{2k+2j+2}{2p})\Gamma(\frac{2k+2j+p+n+4}{2p})}\times\prod\limits_{j=m}^{p+l-1}
\dfrac{(2k+2j+2)\Gamma(\frac{2k+2j+2}{2s})\Gamma(\frac{2k+2j+s+d+4}{2s})}{\Gamma(\frac{2k+2j+4}{2s})\Gamma(\frac{2k+2j+s+d+2}{2s})} \\
&\equiv\frac{\prod\limits_{j=m}^{p+l-1}(2k+2j+2)}{\prod\limits_{j=l}^{m+s-1}(2k+2j+2)} 
\dfrac{\Gamma(\frac{2k+2m+2s+2}{2p})\Gamma(\frac{2k+2l+p+n+2}{2p})}{\Gamma(\frac{2k+2l+2}{2p})\Gamma(\frac{2k+2m+2s+p+n+2}{2p})} \times
\dfrac{\Gamma(\frac{2k+2m+2}{2s})\Gamma(\frac{2k+2p+2l+s+d+2}{2s})}{\Gamma(\frac{2k+2p+2l+2}{2s})\Gamma(\frac{2k+2m+s+d+2}{2s})}   \\
&\equiv\prod\limits_{j=m}^{l-1}(2k+2j+2)
\dfrac{\Gamma(\frac{2k+2l+2p+2}{2p})\Gamma(\frac{2k+2l+p+n+2}{2p})}{\Gamma(\frac{2k+2l+2}{2p})\Gamma(\frac{2k+2l+2p+p+n+2}{2p})} \times
\dfrac{\Gamma(\frac{2k+2m+2}{2s})\Gamma(\frac{2k+2s+2m+s+d+2}{2s})}{\Gamma(\frac{2k+2s+2m+2}{2s})\Gamma(\frac{2k+2m+s+d+2}{2s})}   \\
&\equiv\prod\limits_{j=m}^{l-1}(2k+2j+2)
\dfrac{\Gamma(\frac{2k+2l+2}{2p}+1)\Gamma(\frac{2k+2l+p+n+2}{2p})}{\Gamma(\frac{2k+2l+2}{2p})\Gamma(\frac{2k+2l+p+n+2}{2p}+1)} \times
\dfrac{\Gamma(\frac{2k+2m+2}{2s})\Gamma(\frac{2k+2m+s+d+2}{2s}+1)}{\Gamma(\frac{2k+2m+2}{2s}+1)\Gamma(\frac{2k+2m+s+d+2}{2s})}   \\
&\equiv\prod\limits_{j=m}^{l-1}(2k+2j+2)\dfrac{(2k+2l+2)(2k+2m+s+d+2)}{(2k+2l+p+n+2)(2k+2m+2)},
\end{align*}

\begin{align*}
\prod\limits_{j=0}^{s-1}\widehat{g}(2k+2j+3) &= \prod\limits_{j=0}^{s-1}\dfrac{\Gamma(\frac{2k+2j+4}{2s})\Gamma(\frac{2k+2j+s+d+2}{2s})}{(2k+2j+2)\Gamma(\frac{2k+2j+2}{2s})\Gamma(\frac{2k+2j+s+d+4}{2s})}  \\
&\equiv \prod\limits_{j=0}^{s-1}\dfrac{\Gamma(\frac{2k+2s+2}{2s})\Gamma(\frac{2k+s+d+2}{2s})}{(2k+2j+2)\Gamma(\frac{2k+2}{2s})\Gamma(\frac{2k+2s+s+d+2}{2s})}  \\
&\equiv \dfrac{1}{(2k+s+d+2)\prod\limits_{j=1}^{s-1}(2k+2j+2)},
\end{align*}
and
\begin{equation*}
\prod\limits_{j=0}^{p-1}\widehat{f}(2k+2j+3) \equiv \dfrac{1}{(2k+p+n+2)\prod\limits_{j=1}^{p-1}(2k+2j+2)},
\end{equation*}
it follows that the left hand side "LS" of (26), when taking $z=2k+2$, can be written as
\begin{equation*}
LS \equiv \prod\limits_{j=m}^{l-1}(z+2j)\dfrac{(z+2l)(z+2m+s+d)}{(z+2l+p+n)(z+2m)}\dfrac{(z+p+n)\prod\limits_{j=1}^{p-1}(z+2j)}{(z+s+d)\prod\limits_{j=1}^{s-1}(z+2j)}.
\end{equation*}
On the right side of (\ref{eq20}) we have that
\begin{align*}
& \prod\limits_{j=p}^{m-1}\widehat{f}(2k+2j+3)-\prod\limits_{j=p}^{m-1}\widehat{g}(2k+2j+3) \\
&=  \prod\limits_{j=p}^{m-1}\dfrac{\Gamma(\frac{2k+2j+4}{2p})\Gamma(\frac{2k+2j+p+n+2}{2p})}{(2k+2j+2)\Gamma(\frac{2k+2j+2}{2p})\Gamma(\frac{2k+2j+p+n+4}{2p})}-\prod\limits_{j=p}^{m-1}\dfrac{\Gamma(\frac{2k+2j+4}{2s})\Gamma(\frac{2k+2j+s+d+2}{2s})}{(2k+2j+2)\Gamma(\frac{2k+2j+2}{2s})\Gamma(\frac{2k+2j+s+d+4}{2s})}
 \\
&=\frac{1}{\prod\limits_{j=p}^{m-1}(2k+2j+2)} \Big[\dfrac{\Gamma(\frac{2k+2m+2}{2p})\Gamma(\frac{2k+2p+p+n+2}{2p})}{\Gamma(\frac{2k+2p+2}{2p})\Gamma(\frac{2k+2m+p+n+2}{2p})}-\dfrac{\Gamma(\frac{2k+2m+2}{2s})\Gamma(\frac{2k+2p+s+d+2}{2s})}{\Gamma(\frac{2k+2p+2}{2s})\Gamma(\frac{2k+2m+s+d+2}{2s})}\Big]\\
&= \frac{1}{\prod\limits_{j=p}^{m-1}(z+2j)} \Big[\dfrac{(z+p+n)\Gamma(\frac{z+2m}{2p})\Gamma(\frac{z+p+n}{2p})}{(z)\Gamma(\frac{z}{2p})\Gamma(\frac{z+2m+p+n}{2p})}-\dfrac{\Gamma(\frac{z+2m}{2s})\Gamma(\frac{z+2p+s+d}{2s})}{\Gamma(\frac{z+2p}{2s})\Gamma(\frac{z+2m+s+d}{2s})}\Big]. 
\end{align*}
Finally, Equation (\ref{eq20}) becomes
\begin{equation}
H(z)F(z+2\alpha)\equiv F(z),
\end{equation}
with
\begin{equation*}
H(z)=\dfrac{(z+p+n)(z+2m+s+d)}{(z+2l+p+n)(z+s+d)}
\end{equation*}
and
\begin{equation*}
F(z)=(z+2m)\Big[\dfrac{(z+p+n)\Gamma(\frac{z+2m}{2p})\Gamma(\frac{z+p+n}{2p})}{(z)\Gamma(\frac{z}{2p})\Gamma(\frac{z+2m+p+n}{2p})}-\dfrac{\Gamma(\frac{z+2m}{2s})\Gamma(\frac{z+2p+s+d}{2s})}{\Gamma(\frac{z+2p}{2s})\Gamma(\frac{z+2m+s+d}{2s})}\Big].
\end{equation*}
 Therefore, and using the same argument as in Case I,  we obtain the desired result.
\end{proof}

\end{document}